\documentclass[12pt]{article}%
\usepackage{amsmath,amssymb,amsfonts,amscd,amsthm}
\usepackage{amsfonts}
\usepackage{amsmath}
\usepackage{graphicx}%
\usepackage{float}
\usepackage[margin=1.5cm]{caption}
\usepackage{subcaption}
\usepackage{epstopdf}
\usepackage{xcolor}
\usepackage{multirow}

\usepackage{tikz}
\usetikzlibrary{arrows,shapes,positioning}
\usetikzlibrary{decorations.markings}
\tikzstyle arrowstyle=[scale=1]
\tikzstyle directed=[postaction={decorate,decoration={markings,
    mark=at position .65 with {\arrow[arrowstyle]{stealth}}}}]
\tikzstyle reverse directed=[postaction={decorate,decoration={markings,
    mark=at position .65 with {\arrowreversed[arrowstyle]{stealth};}}}]
    
\providecommand{\U}[1]{\protect\rule{.1in}{.1in}}
\newtheorem{theorem}{Theorem}

\newtheorem{cor}[theorem]{Corollary}

\newtheorem{definition}[theorem]{Definition}

\newtheorem{lemma}[theorem]{Lemma}

\newtheorem{prop}[theorem]{Proposition}

\newcommand{\suchthat}{\;\ifnum\currentgrouptype=16 \middle\fi|\;}

\begin{document}
  \begin{center}
        {\fontsize{18}{22}\selectfont
       \bf Relative equilibria for the positive curved $n$--body problem}
       \end{center}

\vspace{4mm}

        \begin{center}
        {\bf Ernesto P\'erez-Chavela$^1$, and Juan Manuel S\'anchez-Cerritos$^2$}\\
\bigskip
$^1$Departamento de Matem\'aticas\\
Instituto Tecnol\'ogico Aut\'onomo de M\'exico, Mexico City, Mexico\\
\bigskip
$^2$Department of Mathematics\\
Sichuan University, Chengdu, People's Republic of China\\
\bigskip
ernesto.perez@itam.mx, sanchezj01@gmail.com
       \end{center}

        

        \abstract{We consider the $n$--body problem defined on surfaces of constant positive curvature. For the $5$ and $7$--body problem in a collinear symmetric configuration we obtain initial positions which lead to relative equilibria. We give explicitly the values of masses in terms of the initial positions. For positions for which relative equilibria exist, there are infinitely many values of 
the masses that generate such solutions. For the $5$ and $7$--body problem, the set of  parameters (masses and positions) leading to relative equilibria  has positive Lebesgue measure. }
        
\section{Introduction}
We consider the generalization of the gravitational $n-$body problem to spaces of positive constant curvature proposed by Diacu, P\'erez-Chavela and Santoprete \cite{Diacu,Diacu2}. The problem has its roots on the ideas about non-Euclidean geometries proposed by Lovachevski and Bolyai in the 19th century \cite{Bolyai,Lovachevski}. For more details about the history of this fascinating problem we refer the interested readers to \cite{Diacu4}.

In this paper we focus on a special type of solutions, the so called relative equilibria. Roughly speaking, they are solutions of the equations of motion where system of particles moves as a rigid body, or in other words where the mutual distances between the particles remain constant along the time. In the classical gravitational Newtonian case these kind of solutions have been deeply analyzed since the Euler and Lagrange times until present days. Relative equilibria in curved spaces (positive and negative curvature) have been also widely studied in recent years, see for instance \cite{Diacu4,Diacu5,Diacu3,Diacu,paper1,naranjo,Martinez,paper,paperh2,Tibboel,Zhu}, where existence, stability and bifurcations, among others properties of different families have been studied.

In the negative curvature  case, recently the authors analyzed the collinear relative equilibria  for 
$5$ and $7$ bodies on negative curved spaces, founding some interesting results \cite{paperh2}. This paper is a natural continuation on this research line.

A collinear relative equilibria in curved spaces is a relative equilibrium where the particles lie at every time $t$ on the same rotating geodesic.
The results presented on this paper are, as mentioned above, are related with symmetric collinear relative equilibria for the $5$ and $7$--body problem on positive curved spaces. In both cases
we tackle the problem about the distribution of the particles in order to get relative equilibria. 
 We show that for positions that generate relative equilibria there exist infinitely many values of masses that lead to those solutions. The set of parameters (positions and masses) has positive Lebesgue measure. 
 
 We have also get conditions for the no existence of relative equilibria, in fact this result is generalized for the general $n-$symmetrical case. 
        
 \section{Statement of the main results}\label{RE}
        
        Consider a surface of constant curvature $1$. In this paper we use the stereographic model $\mathbb{M}^2$,  given by the complex plane $\mathbb{C}$ endowed with the metric 
        
        \begin{equation}
        ds=\dfrac{4dzd\bar{z}}{(1+|z|^2)^2}, \ \ z=u+iv \in \mathbb{C}.
        \end{equation}

        We denote by $z_i$ the position of the particle with mass $m_i$. The distance  between any two points  in this space  satisfies
        
        \[
        \cos(d(z_k,z_j))=\dfrac{2(z_k \bar{z}_j+z_j\bar{z}_k)+(|z_k|^2-1)(|z_j|-1)}{(|z_k|+1)(|z_j|+1)}.
        \]
        
The potential is given by
\[ U(q)=\sum_{i<j}m_im_j \cot(d(q_i,q_j)). \]

And the kinetic energy is defined by

\[ T=\dfrac{1}{2}\sum_{i}m_i.\dfrac{4}{(1+|z_i|^2)^2}|\dot{z}_i|^2. \]

From the Euler-Lagrange equations, the equations of motion take the form        

\begin{equation}
\ddot{z}_i=\dfrac{2z_i\dot{z_k}^2}{1+|z_i|^2}+\dfrac{(1+|z_i|^2)^2}{2}\dfrac{\partial U}{\partial \bar{z}_i},\label{eqs}
\end{equation}
where

  \[\dfrac{\partial U}{\partial \bar{z}}=\sum_{j=1, j\neq k}^n\dfrac{2m_km_j(1+|z_k|^2)(1+z_j)^2(1+\bar{z}_jz_k)(z_j-z_k)}{\left((|z_k|^2+1)^2(|z_j|^2+1)^2-[2(z_k\bar{z}_j+z_j\bar{z}_k)+(|z_k|^2-1)(|z_j|^2)-1]^2\right)^{3/2}}.\]

The main results of this paper are the following

\begin{theorem}\label{Theorem 1}
In the $5$--body problem on $\mathbb{M}^2$ we consider $5$ particles on the same geodesic with masses $m_1= \mu, m_2=m_3=1$ and $m_4=m_5=m$, and initial positions $z_1=0, \ z_2=-z_3=a>0, \ z_4=-z_5=r>a$ (Figure \ref{fig:cinco}), 
\begin{itemize}
\item If  $ar-1<0$, then do not exist relative equilibria.
\item For any other parameters, there exist relative equilibria.
\end{itemize}
\end{theorem} 

\begin{theorem}\label{Theorem 2}
In the $7$--body problem on $\mathbb{M}^2$, we consider particles on the same geodesic with masses $m_1=\mu, \ m_2=m_3=1, \ m_4=m_5=M. \ m_6=m_7=m$, and initial positions $z_1=0, \ z_2=-z_3=x>0. \ z_4=-z_5=y>0. \ z_6=-z_7=z>0$, ($x<y<z$). The equator of $\mathbb{S}^2$ under the stereographic projection goes to the unit circle on $\mathbb{M}^2$, we call it the geodesic circle.
\begin{itemize}
\item If $m_2,m_3, m_4,m_5,m_6,m_7$ lie inside the geodesic circle, then do not exist relative equilibria.
\item If   $m_2,m_3, m_4,m_5,m_6,m_7$ lie outside the geodesic circle, then there exist initial positions that generate relative equilibria.
\item If $m_2,m_3$ lie inside the geodesic circle, and $ m_4,m_5,m_6,m_7$ lie outside the geodesic circle with $y<z<\frac{1}{x}$, then do not exist relative equilibria. If $y<1/x<z$ or $1/x<y<z$ then it is possible to find relative equilibria.
\item If $m_2,m_3,m_4,m_5$ lie inside the geodesic circle, and $m_6,m_7$ lie outside the geodesic circle with $z<1/y$, then do not exist relative equilibria. If $1/y<z<1/x$ or $1/x<z$ then it is possible to find relative equilibria.
\end{itemize}
\end{theorem}

\bigskip

Before to proceed with the proof of the above results we must point the formal definition of relative equilibria and the frame work in which we will be working.

\medskip

Let $Iso(\mathbb{M}^2)$ be the group of isometries of $\mathbb{M}^2$, and let $\{G(t)  \}$ be a one-parametric subgroup of $Iso(\mathbb{M}^2)$. 

\begin{definition} A Relative Equilibrium  of the curved n-body problem is a solution of (\ref{eqs}) which is invariant relative to the subgroup $\{G(t)  \}$. \end{definition}

It is well known that, in order to obtain relative equilibria on $\mathbb{M}^2$, it is enough to study solutions given by the action $w(t)=e^{it}z(t)$ (see \cite{EPC}), i.e. solutions of the equations of motion where the orbits of the bodies are Euclidean circles. In the same paper, the authors show the necessary conditions for the existence of relative equilibria, which are given by:.

\begin{prop}
Consider $n$ point particles with masses $m_1,m_2, \dots , m_n$ moving on $\mathbb{M}^2$. A necessary and sufficient condition for the solution $z_1,z_2,\dots,z_n$ of (\ref{eqs}) to be a relative equilibrium is that the coordinates  satisfy  the following system given by the rational functions:
\begin{equation}
\dfrac{(1-r_i^2)z_i}{4(1+r_i^2)^4}=-\sum_{j=1,j\neq i}^n \dfrac{m_j(r_j^2+1)^2(1+z_i\bar{z}_j)(z_j-z_i)}{T_{ij}^{3/2}},\label{cond}
\end{equation}
 where
 $T_{ij}=(r_i^2+1)^2(r_j^2+1)^2-[2(z_i\bar{z_j}+z_j\bar{z}_i)+(r_i^2-1)(r_j^2-1)]^2$, and $|z_l|=r_l \in [0,\pi)$.
\end{prop}

It is clear, as they show, that for three particles placed on a geodesic, there exist masses and positions that satisfy such condition. For $n \geq 4$ nothing is know until now. The ideas to extend the results about the existence of relative equilibria for $n \geq 4$, are simple and clear, nevertheless the computations are not easy, as we will show in this work. 
        
\section{Proof of Theorem \ref{Theorem 1}}\label{sec3}  

\begin{figure}[!h]
	\centering
  \begin{tikzpicture}[thick, scale=1.4]
  \draw[-] (0,0) -- (180:1.6);
   \draw[-] (0,0) -- (90:1.6);
   \fill (0,0)  circle[radius=2pt];
   \fill (0.86,0.5)  circle[radius=2pt];
   \fill (-1.6*0.86,-1.6*0.5)  circle[radius=2pt];
   \fill (1.6*0.86,1.6*0.5)  circle[radius=2pt];
   \fill (-0.86,-0.5)  circle[radius=2pt];
   \draw (1,-0.08) -- (1,0.08);
    \draw (-1,-0.08) -- (-1,0.08);
     \draw (1.6,-.08) -- (1.6,.08);
      \draw (-1.6,-.08) -- (-1.6,.08);
  \draw (0,0) circle (1.4);
\draw [->,red,thick,domain=0:26] plot ({cos(\x)}, {sin(\x)});  
\draw [->,blue,thick,domain=0:28] plot ({1.6*cos(\x)}, {1.6*sin(\x)});  
\draw [->,red,thick,domain=180:206] plot ({cos(\x)}, {sin(\x)});  
\draw [->,blue,thick,domain=180:208] plot ({1.6*cos(\x)}, {1.6*sin(\x)});  
    
  \node[above] at(0.1:0.2) {{\small $\mathbf{q}_1$}};
  \node[above] at (30:1.7) {{\small $\mathbf{q}_4$}};
  \node[above] at (225:1.7) {{\small $\mathbf{q}_5$}};
  \node[above] at (30:1) {{\small $\mathbf{q}_2$}};
  \node[above] at (233:1) {{\small $\mathbf{q}_3$}};
  \node[below] at  (1,0){{\tiny $a$}};
  \node[above] at  (-1,0){{\tiny $-a$}};
  \node[below] at  (1.6,0){{\tiny $r$}};
  \node[above] at  (-1.6,0){{\tiny $-r$}};
  \node[below] at  (1.3,0){{\tiny $1$}};
  
   \draw [<->] (0, -2.2)--(0, 2.2) node (zaxis) [left] {$i$};
         \draw [<->](-2.2, 0)--(2.2, 0) node (yaxis) [right] {$\mathbb{R}$};
  \end{tikzpicture}
  \caption{Five bodies on a geodesic in $\mathbb{M}^2$, with $a<1<r$ at time $t>0$.}
  \label{fig:cinco}
  \end{figure}
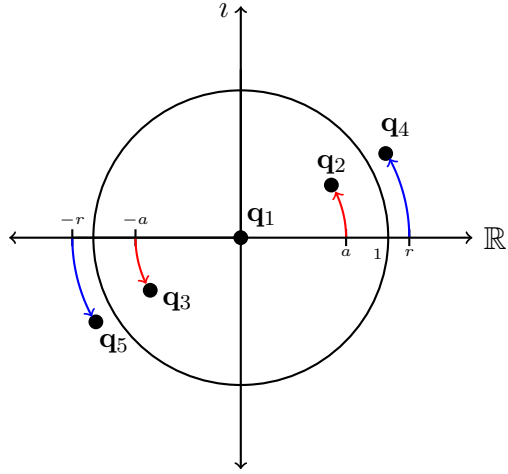

For the symmetric collinear $5$--body problem,  after a suitable rotation, without loss of generality we can consider the initial positions for the configuration as  $z_1=0, \ z_2=-z_3=a, \ z_4=-z_5=r$ (see Figure \ref{fig:cinco}).

 First case. This case corresponds to the particles lying inside the geodesic circle defined in Theorem \ref{Theorem 2}, i.e. $0<a<r<1$ . Using equation (\ref{cond}) for particle $z_2$ and $z_4$, we have

\begin{equation}\label{e1}
\begin{split}
&\dfrac{1}{8}\dfrac{(a^2+1)^2}{a^2(1-a^2)^2}-\dfrac{(1-a^2)a}{(1+a^2)^4}=-\dfrac{1}{2}\dfrac{\mu}{a^2}+\dfrac{mra(r^2+1)(1-r^2)(1-a^2)}{(a^2r^2-1)^2(a^2-r^2)}, \\
&\dfrac{1}{2}(a^2+1)^2\left[ \dfrac{1}{(ar+1)^2(a-r)^2}+\dfrac{1}{(r+a)^2(ar-1)^2} \right]-\dfrac{(1-r^2)r}{(1+r^2)^4}=-\dfrac{1}{2}\dfrac{\mu}{r^2}-\dfrac{1}{8}\dfrac{m(r^2+1)^2}{r(r^2-1)^2}.
\end{split}
\end{equation}

\begin{lemma}\label{lema2}
$\dfrac{1}{2}(a^2+1)^2\left[ \dfrac{1}{(ar+1)^2(a-r)^2}+\dfrac{1}{(r+a)^2(ar-1)^2} \right]-\dfrac{(1-r^2)r}{(1+r^2)^4}>0$.
\end{lemma}

\begin{proof}
Let us define $H= \dfrac{1}{2}\dfrac{(a^2+1)}{(ar-1)^2(a+r)^2}-\dfrac{(1-r^2)r}{(1+r^2)^4}$. We have

\begin{equation}
\dfrac{\partial H}{\partial a}=-\dfrac{(a^2+1)(ar+a+r-1)(ar-a-r-1)}{(r+a)^3(1-ar)^3}.
\end{equation}
Then $\frac{\partial H}{\partial a}=0$ iff $a_1=\dfrac{1-r}{r+1}$ or $a_2=\dfrac{r+1}{r-1}$. We take $a_1$ since $a \in (0,1)$.

The value $a_1$ is a minimum (for a fixed $r$) since $\dfrac{\partial^2 H}{\partial a^2}=\dfrac{4(r+1)^4}{(r^2+1)^4}>0$. Then we compute

\[ H\left(a=\frac{1-r}{a+r},r\right)=\dfrac{(r^2+r+2)(2r^2-r+1)}{(r^2+1)^4}>0, \ \ r=(0,1) .\]

Hence

\[\dfrac{1}{2}(a^2+1)^2\left[ \dfrac{1}{(ar+1)^2(a-r)^2}+\dfrac{1}{(r+a)^2(ar-1)^2} \right]-\dfrac{(1-r^2)r}{(1+r^2)^4}>H>0.\]
\end{proof}

Lemma (\ref{lema2}) implies that second equation of system (\ref{e1}) is never satisfied, since left part is positive and right part is negative for $ a,r<1$.

Second case. We have $1<a<r$. Using condition (\ref{cond}) then

\begin{equation}\label{ss}
\begin{split}
\dfrac{(a^2-1)a}{(1+a^2)^4}&=-\dfrac{1}{2}\dfrac{\mu}{a^2}+\dfrac{1}{8}\dfrac{(a^2+1)^2}{a^2(a^2-1)^2}+\dfrac{1}{2}\dfrac{m(r^2+1)^2}{(ar+1)^2(a-r)^2}+\dfrac{1}{2}\dfrac{m(r^2+1)^2}{(r+a)^2(ar-1)^2},\\
\dfrac{(r^2-1)r}{(1+r^2)^4}&=-\dfrac{1}{2}\dfrac{\mu}{r^2}-\dfrac{1}{2}\dfrac{(a^2+1)^2}{(ar+1)^2(a-r)^2}+\dfrac{1}{2}\dfrac{(a^2+1)^2}{(r+a)^2(ar-1)^2}+\dfrac{1}{8}\dfrac{m(r^2+1)^2}{r^2(r^2-1)^2}.
\end{split}
\end{equation}

We need to see whether or not there exist parameters $a,r,\mu$ and $m$ such that last system is satisfied. Adding both equations of system (\ref{ss}) we have

\begin{equation}\label{ABCD}
A=B\mu+Cm,
\end{equation}
with

\begin{equation}\label{AA}
\begin{split}
A&=\dfrac{(a^2-1)a}{(a^2+1)^4}+\dfrac{(1-r^2)r}{(1+r^2)^4}\\
&+\dfrac{1}{2}(a^2+1)^2\left( -\dfrac{1}{4}\dfrac{1}{a^2(a^2-1)^2}+\dfrac{1}{(ar+1)^2(a-r)^2}-\dfrac{1}{(r+a)^2(ar-1)}\right),\\
B&=-\dfrac{1}{2}\left( \dfrac{1}{a^2}+\dfrac{1}{r^2}  \right)<0,\\
C&=\dfrac{1}{2}(r^2+1)^2\left(\dfrac{1}{(ar+1)^2(a-r)^2}+\dfrac{1}{(r+a)^2(ar-1)^2}+\dfrac{1}{4}\dfrac{1}{r^2(r^2-1)^2}   \right)>0.
\end{split}
\end{equation}

The sign of $A$ can be positive or negative. From equation (\ref{ABCD}) we have 
\begin{equation}\label{A1}
A-B\mu=Cm.
\end{equation}

We have that there exist masses that generate relative equilibria if $A-B\mu>0$ or

\[ \dfrac{A}{B}<\mu. \]

The mass relation is given by

\[ m=\dfrac{A-B\mu}{C}. \]

Third case. We have $0<a<1<r$. 
%

Here we have two sub cases, when $ar-1<0$ or $ar-1>0$. Recall that $ar-1=0$ correspond to a singularity of the equations of motion.

Consider first the subcase $ar-1<0$.

\begin{equation}\label{sy0}
\begin{split}
-\dfrac{(1-a^2)a}{(1+a^2)^4}&=-\dfrac{1}{2}\dfrac{\mu}{a^2}-\dfrac{1}{8}\dfrac{(a^2+1)^2}{a^2(1-a^2)^2}-\dfrac{2m(r^2+1)^2ar(r^2-1)(1-a^2)}{(a^2r^2-1)^2(r^2-a^2)^2},\\
\dfrac{(r^2-1)r}{(1+r^2)^4}&=-\dfrac{1}{2}\dfrac{\mu}{r^2}-\dfrac{1}{2}\dfrac{(a^2+1)^2}{(ar+1)^2(a-r)^2}-\dfrac{1}{2}\dfrac{(a^2+1)^2}{(r+a)^2(ar-1)^2}+\dfrac{1}{8}\dfrac{m(r^2+1)^2}{r^2(r^2-1)^2}.
\end{split}
\end{equation}

\begin{lemma}\label{lemma5}
$-\dfrac{(1-x^2)x}{(1+x^2)^4}+1/8\,{
\frac { \left( {x}^{2}+1 \right) ^{2}}{ \left( {x}^{2}-1 \right) ^{2}{
x}^{2}}}>0$, $x\in(0,1)$
\end{lemma}

\begin{proof}
We have

\[ -\dfrac{(1-x^2)x}{(1+x^2)^4}+1/8\,{
\frac { \left( {x}^{2}+1 \right) ^{2}}{ \left( {x}^{2}-1 \right) ^{2}{
x}^{2}}}=\dfrac{1}{8}\dfrac{f(x)g(x)}{x^2(x^2-1)^2(x^2+1)^4}, \]
with

$f(x)=x^4+2x^3+2x^2-2x+1>0$, $g(x)=h(x)+D(x)$ where $h(x)=x^6(x^2-2x+8)$ and $D(x)=-2x^5-2x^4+2x^3+8x^2+2x+1$. The function $D(x)$ has only one positive root (Descartes' rule of signs) between $x=1$ and $x=2$, we also have $D(0)=1$. The function $h(x)$ is easy to check that is positive if $x>0$.

Since $f(x)$ is also positive, we conclude $A_1>0$.

\end{proof}

By the above lemma, for this subcase,  the first equation of system (\ref{sy0}) has no solution.
Hence there are not relative equilibria for this positions.

\medskip

Now consider the subcase $ar-1>0$

\begin{equation}\label{sy1}
\begin{split}
-\dfrac{(1-a^2)a}{(1+a^2)^4}&+\dfrac{1}{8}\dfrac{(a^2+1)^2}{a^2(1-a^2)^2}=-\dfrac{1}{2}\dfrac{\mu}{a^2}\\
&+\dfrac{m(r^2+1)^2}{2}\left(  \dfrac{1}{(ar+1)^2(a-r)^2}+\dfrac{1}{(r+a)^2(ar-1)^2} \right),\\
\dfrac{(r^2-1)r}{(1+r^2)^4}&-\dfrac{2ar(a^2+1)^2(r^2-1)(1-a^2)}{(a^2r^2-1)^2(a^2-r^2)^2}=-\dfrac{1}{2}\dfrac{\mu}{r^2}+\dfrac{1}{8}\dfrac{m(r^2+1)^2}{r^2(r^2-1)^2}.\\
\end{split}
\end{equation}

 If we add both equations of system (\ref{sy1}), then we have

\begin{equation}\label{F}
F_1=F_2\mu+F_3m,
\end{equation}
with
\begin{equation}\label{FF}
\begin{split}
F_1&=-\dfrac{(1-a^2)a}{(1+a^2)^4}+\dfrac{1}{8}\dfrac{(a^2+1)^2}{a^2(1-a^2)^2}+\dfrac{(r^2-1)r}{(1+r^2)^4}-\dfrac{2ar(a^2+1)^2(r^2-1)(1-a^2)}{(a^2r^2-1)^2(a^2-r^2)^2}\\
F_2&=-\dfrac{1}{2}\left( \dfrac{1}{a^2}+\dfrac{1}{r^2}  \right)<0,\\
F_3&=\dfrac{(r^2+1)^2}{2}\left(  \dfrac{1}{(ar+1)^2(a-r)^2}+\dfrac{1}{(r+a)^2(ar-1)^2} \right)+\dfrac{1}{8}\dfrac{(r^2+1)^2}{r^2(r^2-1)^2}>0.\\
\end{split}
\end{equation}

The sign of $F_1$ can be positive or negative. From equation (\ref{F}),

\[ F_3m=F_1-F_2\mu. \]

We have that there exist relative equilibria if $F_1-F_2\mu>0$ or 

\[ \dfrac{F_1}{F_2}<\mu .\]

The mass relation is 

\[  m=\dfrac{F_1-F_2\mu}{F_3}.\]

We summarize the conditions in the following table. Consider the values in (\ref{AA}) and (\ref{FF}).
\renewcommand{\arraystretch}{2.9}
\begin{center}
\begin{tabular}{ |c|c|c|c| } 
\hline
 & \multicolumn{2}{|c|}{Positions} & Masses  \\
\hline
$a<r<1$ & \multicolumn{3}{|c|}{No relative equilibria} \\
\hline
\multirow{3}{5em}{$1<a<r$} & \multicolumn{2}{|c|}{$A\geq 0$} & $\mu \in \mathbb{R}^+$, \ \  $m=\dfrac{A-B\mu}{C}$ \\ \cline{2-4} 
						& \multicolumn{2}{|c|}{$A<0$} & $\dfrac{A}{B}<\mu$, \ \  $m=\dfrac{A-B\mu}{C}$ \\ \cline{2-3}
\hline
\multirow{4}{5em}{$a<1<r$} & \multicolumn{2}{|c|}{$ar-1<0$} & No relative equilibria  \\ \cline{2-4}
						& $ar-1>0$ &  $F_1\geq0$  & $\mu \in \mathbb{R}^+$, \ \ $m=\dfrac{F_1-F_2\mu}{F_3}$ \\  \cline{3-4}
						&          &  $F_1<0$ &   $\dfrac{F_1}{F_2}<\mu$, \ \ $m=\dfrac{F_1-F_2\mu}{F_3}$ \\  \cline{3-4}
\hline
\end{tabular}
\end{center}

\medskip

\begin{cor}\label{4b}
In the 4-body problem on $\mathbb{M}^2$ we consider 4 particles on the same geodesic with masses $m_1=m_2=1$ and $m_3=m_4=m$, in a symmetric configuration with initial positions $z_1=-z_2=a>0$ and $z_3=-z_4=r>a$. Then do not exist relative equilibria.
\end{cor}

\begin{proof}
It is enough to analyze the cases $a<r<1$ and $a<1<r$ with $ar<1$. The proof is similar as in the previous theorem, by considering $\mu=0$. It is enough to analyze the case $a<r<1$ and if $a<1<r$ the case $ar-1<0$.

\begin{itemize}
\item Case $a<r<1$
\end{itemize}
The condition (\ref{cond}) is 
\begin{equation}\label{e1}
\begin{split}
&\dfrac{1}{8}\dfrac{(a^2+1)^2}{a^2(1-a^2)^2}-\dfrac{(1-a^2)a}{(1+a^2)^4}=\dfrac{mra(r^2+1)(1-r^2)(1-a^2)}{(a^2r^2-1)^2(a^2-r^2)}, \\
&\dfrac{1}{2}(a^2+1)^2\left[ \dfrac{1}{(ar+1)^2(a-r)^2}+\dfrac{1}{(r+a)^2(ar-1)^2} \right]-\dfrac{(1-r^2)r}{(1+r^2)^4}=-\dfrac{1}{8}\dfrac{m(r^2+1)^2}{r(r^2-1)^2}.
\end{split}
\end{equation}
Lemma \ref{lema2} implies that there is no solution for the second equation.

\begin{itemize}
\item Case $a<1<r$, $ar<1$
\end{itemize}

Condition \ref{cond} is
\begin{equation}
\begin{split}
-\dfrac{(1-a^2)a}{(1+a^2)^4}&=-\dfrac{1}{8}\dfrac{(a^2+1)^2}{a^2(1-a^2)^2}-\dfrac{2m(r^2+1)^2ar(r^2-1)(1-a^2)}{(a^2r^2-1)^2(r^2-a^2)^2},\\
\dfrac{(r^2-1)r}{(1+r^2)^4}&=-\dfrac{1}{2}\dfrac{(a^2+1)^2}{(ar+1)^2(a-r)^2}-\dfrac{1}{2}\dfrac{(a^2+1)^2}{(r+a)^2(ar-1)^2}+\dfrac{1}{8}\dfrac{m(r^2+1)^2}{r^2(r^2-1)^2}.
\end{split}
\end{equation}

We previously checked that 
 
\[ -\dfrac{(1-a^2)a}{(1+a^2)^4}+\dfrac{1}{8}\dfrac{(a^2+1)^2}{a^2(1-a^2)^2}>0. \]

 Hence first equation has no solutions and we conclude that there are not relative equilibria.
\end{proof}

\section{Proof of Theorem \ref{Theorem 2}}\label{sec4} 
\begin{itemize}
\item Case: $m_2,m_3, m_4,m_5,m_6,m_7$ lie inside the geodesic circle. 
\end{itemize}

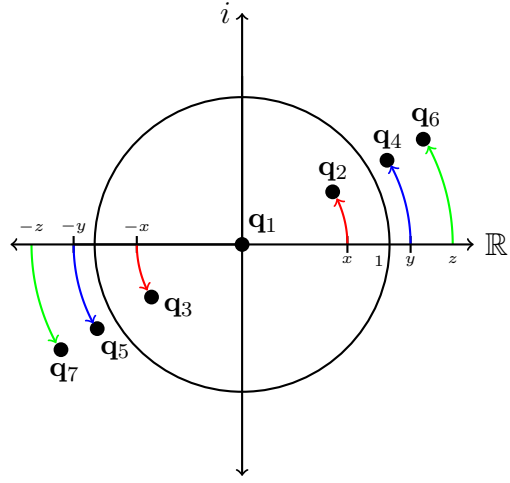
\begin{figure}[!h]
	\centering
  \begin{tikzpicture}[thick, scale=1.4]
  \draw[-] (0,0) -- (180:1.6);
   \draw[-] (0,0) -- (90:1.6);
   \fill (0,0)  circle[radius=2pt];
   \fill (0.86,0.5)  circle[radius=2pt];
   \fill (-1.6*0.86,-1.6*0.5)  circle[radius=2pt];
   \fill (1.6*0.86,1.6*0.5)  circle[radius=2pt];
   \fill (-0.86,-0.5)  circle[radius=2pt];
   \fill (2*0.86,2*0.5)  circle[radius=2pt];
   \fill (-2*0.86,-2*0.5)  circle[radius=2pt];
   \draw (1,-0.08) -- (1,0.08);
    \draw (-1,-0.08) -- (-1,0.08);
     \draw (1.6,-.08) -- (1.6,.08);
      \draw (-1.6,-.08) -- (-1.6,.08);
  \draw (0,0) circle (1.4);
\draw [->,red,thick,domain=0:26] plot ({cos(\x)}, {sin(\x)});  
\draw [->,blue,thick,domain=0:28] plot ({1.6*cos(\x)}, {1.6*sin(\x)});  
\draw [->,red,thick,domain=180:206] plot ({cos(\x)}, {sin(\x)});  
\draw [->,blue,thick,domain=180:208] plot ({1.6*cos(\x)}, {1.6*sin(\x)});  
\draw [->,green,thick,domain=0:28] plot ({2*cos(\x)}, {2*sin(\x)}); 
\draw [->,green,thick,domain=180:208] plot ({2*cos(\x)}, {2*sin(\x)}); 
    
  \node[above] at(0.1:0.2) {{\small $\mathbf{q}_1$}};
  \node[above] at (30:1.6) {{\small $\mathbf{q}_4$}};
  \node[above] at (225:1.7) {{\small $\mathbf{q}_5$}};
  \node[above] at (30:1) {{\small $\mathbf{q}_2$}};
  \node[above] at (233:1) {{\small $\mathbf{q}_3$}};
  \node[above] at (30:2.02) {{\small $\mathbf{q}_6$}};
  \node[above] at (220:2.2) {{\small $\mathbf{q}_7$}};
  \node[below] at  (1,0){{\tiny $x$}};
  \node[above] at  (-1,0){{\tiny $-x$}};
  \node[below] at  (1.6,0){{\tiny $y$}};
  \node[above] at  (-1.6,0){{\tiny $-y$}};
  \node[below] at  (2,0){{\tiny $z$}};
  \node[above] at  (-2,0){{\tiny $-z$}};
  \node[below] at  (1.3,0){{\tiny $1$}};
  
   \draw [<->] (0, -2.2)--(0, 2.2) node (zaxis) [left] {$i$};
         \draw [<->](-2.2, 0)--(2.2, 0) node (yaxis) [right] {$\mathbb{R}$};
  \end{tikzpicture}
  \caption{Seven bodies on a geodesic in $\mathbb{M}^2$, with $x<1<y<z$ at time $t>0$.}
  \label{fig:siete}
  \end{figure}

We start by  considering $m_2,m_3, m_4,m_5,m_6,m_7$ inside the geodesic circle. This case correspond to $x<y<z<1$. Using condition (\ref{cond}) for particle $z_6$ we obtain (see Figure \ref{fig:siete}).

\begin{equation}
\begin{split}
-\dfrac{(1-z^2)z}{(1+z^2)^4} & +\dfrac{1}{2}(x^2+1)^2\left(\dfrac{1}{(xz+1)^2(x-z)^2}+\dfrac{1}{(z+x)^2(xz-1)^2}  \right)=-\dfrac{1}{2}\dfrac{\mu}{z^2}\\
&-\dfrac{1}{2}(y^2+1)^2\left( \dfrac{1}{(yz+1)^2}+\dfrac{1}{(z+y)^2(yz-1)^2} \right)M-\dfrac{1}{8}\dfrac{(z^2+1)^2m}{(z^2-1)^2z^2}.
\end{split}
\end{equation}

Last equation is never satisfied for $\mu,M,m>0$, since left part is positive. The proof where left part is positive is similar as in lemma (\ref{lema2}).

\begin{itemize}
\item Case: $m_2,m_3, m_4,m_5,m_6,m_7$ lie outside the geodesic circle. 
\end{itemize}

In this case, system (\ref{cond}) takes the form

\begin{equation}\label{st}
\begin{split}
\dfrac{(x^2-1)x}{(1+x^2)^4}=&-1/2\,{\frac {\mu}{{x}^{2}}}+ \left( 1/2\,{\frac { \left( {y}^{2}
+1 \right) ^{2}}{ \left( xy+1 \right) ^{2} \left( -y+x \right) ^{2}}}+
1/2\,{\frac { \left( {y}^{2}+1 \right) ^{2}}{ \left( y+x \right) ^{2}
 \left( xy-1 \right) ^{2}}} \right) {M}\\
 &+ \left( 1/2\,{\frac {
 \left( {z}^{2}+1 \right) ^{2}}{ \left( xz+1 \right) ^{2} \left( x-z
 \right) ^{2}}}+1/2\,{\frac { \left( {z}^{2}+1 \right) ^{2}}{ \left( z
+x \right) ^{2} \left( xz-1 \right) ^{2}}} \right) {m}\\
&+1/8\,{
\frac { \left( {x}^{2}+1 \right) ^{2}}{ \left( {x}^{2}-1 \right) ^{2}{
x}^{2}}},\\
\dfrac{(y^2-1)y}{(1+y^2)^4}=&-1/2\,{\frac {{\mu}}{{y}^{2}}}+1/8\,{\frac {{M}\, \left( {y}^
{2}+1 \right) ^{2}}{ \left( {y}^{2}-1 \right) ^{2}{y}^{2}}}\\
&+ \left( 1/
2\,{\frac { \left( {z}^{2}+1 \right) ^{2}}{ \left( yz+1 \right) ^{2}
 \left( -z+y \right) ^{2}}}+1/2\,{\frac { \left( {z}^{2}+1 \right) ^{2
}}{ \left( z+y \right) ^{2} \left( yz-1 \right) ^{2}}} \right) {m
}\\
&-1/2\,{\frac { \left( {x}^{2}+1 \right) ^{2}}{ \left( xy+1 \right) ^{
2} \left( -y+x \right) ^{2}}}+1/2\,{\frac { \left( {x}^{2}+1 \right) ^
{2}}{ \left( y+x \right) ^{2} \left( xy-1 \right) ^{2}}},\\
\dfrac{(z^2-1)z}{(1+z^2)^4}=&-1/2\,{\frac {\mu}{{z}^{2}}}+ \left( -1/2\,{\frac { \left( {y}^{2
}+1 \right) ^{2}}{ \left( yz+1 \right) ^{2} \left( -z+y \right) ^{2}}}
+1/2\,{\frac { \left( {y}^{2}+1 \right) ^{2}}{ \left( z+y \right) ^{2}
 \left( yz-1 \right) ^{2}}} \right) {M}\\
 &-1/2\,{\frac { \left( {x}^
{2}+1 \right) ^{2}}{ \left( xz+1 \right) ^{2} \left( x-z \right) ^{2}}
}+1/2\,{\frac { \left( {x}^{2}+1 \right) ^{2}}{ \left( z+x \right) ^{2
} \left( xz-1 \right) ^{2}}}+1/8\,{\frac {{m}\, \left( {z}^{2}+1
 \right) ^{2}}{ \left( {z}^{2}-1 \right) ^{2}{z}^{2}}}.
\end{split}
\end{equation}

Let be

\begin{equation}
\begin{split}
A_1=&\dfrac{(x^2-1)x}{(1+x^2)^4}-1/8\,{
\frac { \left( {x}^{2}+1 \right) ^{2}}{ \left( {x}^{2}-1 \right) ^{2}{
x}^{2}}},\\
A_2=&\dfrac{(y^2-1)y}{(1+y^2)^4}+1/2\,{\frac { \left( {x}^{2}+1 \right) ^{2}}{ \left( xy+1 \right) ^{
2} \left( -y+x \right) ^{2}}}-1/2\,{\frac { \left( {x}^{2}+1 \right) ^
{2}}{ \left( y+x \right) ^{2} \left( xy-1 \right) ^{2}}},\\
A_3=&\dfrac{(z^2-1)z}{(1+z^2)^4}+1/2\,{\frac { \left( {x}^
{2}+1 \right) ^{2}}{ \left( xz+1 \right) ^{2} \left( x-z \right) ^{2}}
}-1/2\,{\frac { \left( {x}^{2}+1 \right) ^{2}}{ \left( z+x \right) ^{2
} \left( xz-1 \right) ^{2}}},\\
a_1=&-1/2\,{\frac {1}{{x}^{2}}},\\
a_2=&-1/2\,{\frac {1}{{y}^{2}}},\\
a_3=&-1/2\,{\frac {1}{{z}^{2}}},\\
b_1=& 1/2\,{\frac { \left( {y}^{2}
+1 \right) ^{2}}{ \left( xy+1 \right) ^{2} \left( -y+x \right) ^{2}}}+
1/2\,{\frac { \left( {y}^{2}+1 \right) ^{2}}{ \left( y+x \right) ^{2}
 \left( xy-1 \right) ^{2}}},\\
 b_2=&1/8\frac { \left( {y}^
{2}+1 \right) ^{2}}{ \left( {y}^{2}-1 \right) ^{2}{y}^{2}},\\
b_3=&-1/2\,{\frac { \left( {y}^{2
}+1 \right) ^{2}}{ \left( yz+1 \right) ^{2} \left( -z+y \right) ^{2}}}
+1/2\,{\frac { \left( {y}^{2}+1 \right) ^{2}}{ \left( z+y \right) ^{2}
 \left( yz-1 \right) ^{2}}},
\end{split}
\end{equation}

\begin{equation}
\begin{split}
c_1=&1/2\,{\frac {
 \left( {z}^{2}+1 \right) ^{2}}{ \left( xz+1 \right) ^{2} \left( x-z
 \right) ^{2}}}+1/2\,{\frac { \left( {z}^{2}+1 \right) ^{2}}{ \left( z
+x \right) ^{2} \left( xz-1 \right) ^{2}}},\\
c_2=&1/
2\,{\frac { \left( {z}^{2}+1 \right) ^{2}}{ \left( yz+1 \right) ^{2}
 \left( -z+y \right) ^{2}}}+1/2\,{\frac { \left( {z}^{2}+1 \right) ^{2
}}{ \left( z+y \right) ^{2} \left( yz-1 \right) ^{2}}},\\
c_3=&1/8\,{\frac {{m}\, \left( {z}^{2}+1
 \right) ^{2}}{ \left( {z}^{2}-1 \right) ^{2}{z}^{2}}}.
\end{split}
\end{equation}

System (\ref{st}) becomes

\begin{equation}
\begin{split}
A_1=a_1\mu+b_1M+c_1m,\\
A_2=a_2\mu+b_2M+c_2m,\\
A_3=a_3\mu+b_3M+c_3m.\\
\end{split}
\end{equation}

\begin{lemma}\label{lemma4}
$\dfrac{(x^2-1)x}{(1+x^2)^4}-1/8\,{
\frac { \left( {x}^{2}+1 \right) ^{2}}{ \left( {x}^{2}-1 \right) ^{2}{
x}^{2}}}<0$, $x>1$.
\end{lemma}

\begin{proof}
The function $A_1$ can be written as 

\[ A_1=-\dfrac{1}{8}\dfrac{f(x)g(x)}{x^2(x^2-1)^2(x^2+1)^4}, \]
with

$f(x)=x^4-2x^3+2x^2+2x+1>0$, $g(x)=x^8+2x^7+8x^6+h(x)$, $h(x)=3xD(x)$ and $D(x)=x^4-x^3-x^2+4x-1$. The function $D(x)$ has a minimum (its only critical point) at $x=-\frac{1}{12}(729+12\sqrt{3441})^{1/3}-\dfrac{11}{4(729+12\sqrt{3441})^{1/3}}+\frac{1}{4}\approx-0.9334$. We have $D(1)=2$, hence $D(x)>2$ for $x>1$. We can conclude that $h(x)$ and $g(x)$ are positive. It is not difficult to check that $f(x)$ is also positive. All this facts implies that $A_1<0$.

\end{proof}

The above lemma implies $A_1<0$. From the first equation we have $\mu=\dfrac{A_1}{a_1}-\dfrac{b_1}{a_1}M -\dfrac{c_1}{a_1}m>0$. Substituting this value into the other equations, then

\begin{equation}
\begin{split}
\dfrac{a_2A_1}{a_1}+\left[  -\dfrac{a_2b_1}{a_1}+b_2 \right]M+\left[ -\dfrac{a_2c_1}{a_1}+c_2 \right]m=A_2,\\
\dfrac{a_3A_1}{a_1}+\left[  -\dfrac{a_3b_1}{a_1}+b_3 \right]M+\left[ -\dfrac{a_3c_1}{a_1}+c_3 \right]m=A_3.
\end{split}
\end{equation}

We need to see if last system has positive solutions for $M$ and $m$.

Adding both equations of last system we get

\begin{equation}\label{eee}
\begin{split}
 A_1[a_2+a_3]-a_1[A_2+A_3]+&M\left( b_1[-a_2-a_3]+a_1[b_2+b_3]\right)\\
&+m\left(c_1[-a_2-a_3]+a_1[c_2+c_3]\right)=0 .
\end{split}
\end{equation}

We have $A_1,a_1,a_2,a_3<0$; $b_1,b_2,b_3,c_1,c_2,c_3>0$. Let us fix $x$ and $z$. Notice that the only functions that depend on $y$ and $z$ simultaneously are $b_3$ and $c_2$. Consider values of $y$ close enough to $z$ in such a way that $b_1[-a_2-a_3]+a_1[b_2+b_3]<0$  and  $c_1[-a_2-a_3]+a_1[c_2+c_3]<0$ (since $a_1>0$). This is enough to conclude that there exist $M$ and $\mu$ such that (\ref{eee}) is satisfied.

\begin{itemize}
\item Case: $m_2, m_3$ lie inside the geodesic circle and $m_4,m_5,m_6,m_7$ lie outside the geodesic circle, with $y<z<1/x$.
\end{itemize}
The equation for this case corresponding to particle $m_2$ is

\begin{equation}\label{r1}
\begin{split}
-\dfrac{(1-x^2)x}{(1+x^2)^4}&+1/8\,{
\frac { \left( {x}^{2}+1 \right) ^{2}}{ \left( {x}^{2}-1 \right) ^{2}{
x}^{2}}}=-1/2\,{\frac {\mu}{{x}^{2}}}\\
&+ \left( 1/2\,{\frac { \left( {y}^{2}
+1 \right) ^{2}}{ \left( xy+1 \right) ^{2} \left( -y+x \right) ^{2}}}-
1/2\,{\frac { \left( {y}^{2}+1 \right) ^{2}}{ \left( y+x \right) ^{2}
 \left( xy-1 \right) ^{2}}} \right) {M}\\
 &+ \left( 1/2\,{\frac {
 \left( {z}^{2}+1 \right) ^{2}}{ \left( xz+1 \right) ^{2} \left( x-z
 \right) ^{2}}}-1/2\,{\frac { \left( {z}^{2}+1 \right) ^{2}}{ \left( z
+x \right) ^{2} \left( xz-1 \right) ^{2}}} \right) {m}.\\
\end{split}
\end{equation}

The factors for $m$ and $M$ in the last equation can be seen as

\[ \dfrac{1}{2}\dfrac{yx(y^2+1)^2(y^2-1)(x^2-1)}{(xy+1)^2(x^2-y^2)^2(xy-1)^2}<0, \]
\[  \dfrac{1}{2}\dfrac{zx(z^2+1)^2(z^2-1)(x^2-1)}{(xz+1)^2(x^2-z^2)^2(xz-1)^2}<0.
\]

Hence (\ref{r1}) is never satisfied, since left part of the equation is positive (Lemma \ref{lemma5}).

\begin{itemize}
\item Case: $m_2, m_3$ lie inside the geodesic circle and$ m_4,m_5,m_6,m_7$ lie outside the geodesic circle, with $y<1/x<z$.
\end{itemize}
The equations of motion become

\begin{equation}\label{ss2}
\begin{split}
-\dfrac{(1-x^2)x}{(1+x^2)^4}+&1/8\,{\frac { \left( {x}^{2}+1 \right) ^{2}}{ \left( {x}^{2}-1
 \right) ^{2}{x}^{2}}}
=-1/2\,{\frac {{\mu}}{{x}^{2}}}- \dfrac{2yx(y^2+1)^2(y^2-1)(1-x^2)}{(x^2y^2-1)^2(x^2-y^2)^2} {M}\\
 &+ \left( 1/2\,{\frac {
 \left( {z}^{2}+1 \right) ^{2}}{ \left( xz+1 \right) ^{2} \left( x-z
 \right) ^{2}}}+1/2\,{\frac { \left( {z}^{2}+1 \right) ^{2}}{ \left( z
+x \right) ^{2} \left( xz-1 \right) ^{2}}} \right) {m},
\\
\dfrac{(y^2-1)y}{(1+y^2)^4}&+1/2\,{\frac { \left( {x}^{2}+1 \right) ^{2}}{ \left( y+x \right) ^{2}
 \left( xy-1 \right) ^{2}}}+1/2\,{\frac { \left( {x}^{2}+1 \right) ^{2
}}{ \left( xy+1 \right) ^{2} \left( x-y \right) ^{2}}}
=-1/2\,{\frac {{\mu}}{{y}^{2}}}\\
&+1/8\,{\frac {{M}\, \left( {y}^
{2}+1 \right) ^{2}}{ \left( {y}^{2}-1 \right) ^{2}{y}^{2}}}+ \left( 1/
2\,{\frac { \left( {z}^{2}+1 \right) ^{2}}{ \left( yz+1 \right) ^{2}
 \left( -z+y \right) ^{2}}}+1/2\,{\frac { \left( {z}^{2}+1 \right) ^{2
}}{ \left( z+y \right) ^{2} \left( yz-1 \right) ^{2}}} \right) {m
},
\\
\dfrac{(z^2-1)z}{(1+z^2)^4}-&\dfrac{2zx(x^2+1)^2(1-x^2)(z^2-1)}{(x^2z^2-1)^2(x^2-z^2)^2}=-1/2\,{\frac {{\mu}}{{z}^{2}}}\\
& -\dfrac{2zy(y^2+1)^2(y^2-1)(z^2-1)}{(y^2z^2-1)^2(y^2-z^2)^2}{M}+1/8\,{\frac {{m}\,
 \left( {z}^{2}+1 \right) ^{2}}{ \left( {z}^{2}-1 \right) ^{2}{z}^{2}}}
\end{split}
\end{equation}

Let be

\begin{equation}
\begin{split}
A_1=&-\dfrac{(1-x^2)x}{(1+x^2)^4}+1/8\,{\frac { \left( {x}^{2}+1 \right) ^{2}}{ \left( {x}^{2}-1
 \right) ^{2}{x}^{2}}},\\
 A_2=&\dfrac{(y^2-1)y}{(1+y^2)^4}+1/2\,{\frac { \left( {x}^{2}+1 \right) ^{2}}{ \left( y+x \right) ^{2}
 \left( xy-1 \right) ^{2}}}+1/2\,{\frac { \left( {x}^{2}+1 \right) ^{2
}}{ \left( xy+1 \right) ^{2} \left( x-y \right) ^{2}}},\\
A_3=&\dfrac{(z^2-1)z}{(1+z^2)^4}-\dfrac{2zx(x^2+1)^2(1-x^2)(z^2-1)}{(x^2z^2-1)^2(x^2-z^2)^2}\\
a_1=&-1/2\,{\frac {1}{{x}^{2}}},\\
a_2=&-1/2\,{\frac {1}{{y}^{2}}},\\
a_3=&-1/2\,{\frac {1}{{z}^{2}}},\\
b_1=&- \dfrac{2yx(y^2+1)^2(y^2-1)(1-x^2)}{(x^2y^2-1)^2(x^2-y^2)^2},\\
b_2=&1/8\,{\frac { \left( {y}^
{2}+1 \right) ^{2}}{ \left( {y}^{2}-1 \right) ^{2}{y}^{2}}},\\
b_3=&-\dfrac{2zy(y^2+1)^2(y^2-1)(z^2-1)}{(y^2z^2-1)^2(y^2-z^2)^2},\\
c_1=&\left( 1/2\,{\frac {
 \left( {z}^{2}+1 \right) ^{2}}{ \left( xz+1 \right) ^{2} \left( x-z
 \right) ^{2}}}+1/2\,{\frac { \left( {z}^{2}+1 \right) ^{2}}{ \left( z
+x \right) ^{2} \left( xz-1 \right) ^{2}}} \right),\\
c_2=&\left( 1/
2\,{\frac { \left( {z}^{2}+1 \right) ^{2}}{ \left( yz+1 \right) ^{2}
 \left( -z+y \right) ^{2}}}+1/2\,{\frac { \left( {z}^{2}+1 \right) ^{2
}}{ \left( z+y \right) ^{2} \left( yz-1 \right) ^{2}}} \right), \\
c_3=&1/8\,{\frac {
 \left( {z}^{2}+1 \right) ^{2}}{ \left( {z}^{2}-1 \right) ^{2}{z}^{2}}}.\\
\end{split}
\end{equation}

We have $a_1,a_2,a_3,b_1,b_3<0$; $A_1, A_2,b_2,c_1,c_2,c_3>0$ ($A_1$ is positive by Lemma \ref{lemma5}). System (\ref{ss2}) becomes

\begin{equation}\label{system1}
\begin{split}
A_1=a_1\mu+b_1M+c_1m,\\
A_2=a_2\mu+b_2M+c_2m,\\
A_3=a_3\mu+b_3M+c_3m.\\
\end{split}
\end{equation}

From the first equation we have $m=\dfrac{A_1}{c_1}-\dfrac{a_1\mu}{c_1}-\dfrac{b_1M}{c_1}>0$. Substituting into the other two equations we have

\begin{equation}
\begin{split}
c_1A_2-c_2A_1=\mu(a_2c_1-c_2a_1)+M(b_2c_1-c_2b_1),\\
c_1A_3-c_3A_1=\mu(a_3c_1-c_3a_1)+M(b_3c_1-c_3b_1).
\end{split}
\end{equation}

Adding last two equations,

\begin{equation}\label{Mm}
\begin{split}
-c_1(A_2+A_3)+A_1(c_2+c_3)&+\mu[c_1(a_2+a_3)-a_1(c_2+c_3)]\\
&+M[c_1(b_2+b_3)-b_1(c_2+c_3)]=0.
\end{split}
\end{equation}

Notice that among the functions $A_i,a_i,b_i,c_i, i=1,2,3$, the only values that depend on $x$ and $y$ simultaneously are $A_2$ and $b_1$. Consider values of $y$ close enough to $\frac{1}{x}$, in such a way that $-c_1(A_2+A_3)+A_1(c_2+c_3)<0$ and $[c_1(b_2+b_3)-b_1(c_2+c_3)]>0$ are satisfied.  When these two inequalities are fulfilled, then we can conclude that there are values for $M$ and $\mu$ such that equation (\ref{Mm}) is valid.

\begin{itemize}
\item Case: $m_2, m_3$ lie inside the geodesic circle and$ m_4,m_5,m_6,m_7$ lie outside the geodesic circle, with $1/x<y<z$.
\end{itemize}

The  condition (\ref{cond}) become

\begin{equation}\label{ss3}
\begin{split}
-\dfrac{(1-x^2)x}{(1+x^2)^4}&+1/8\,{\frac { \left( {x}^{2}+1 \right) ^{2}}{ \left( {x}^{2}-1
 \right) ^{2}{x}^{2}}}=-1/2\,{\frac {{\mu}}{{x}^{2}}}\\
& + \left( 1/2\,{\frac { \left( {y}^{2}
+1 \right) ^{2}}{ \left( y\,x+1 \right) ^{2} \left( x-y \right) ^{2}}}
+1/2\,{\frac { \left( {y}^{2}+1 \right) ^{2}}{ \left( y+x \right) ^{2}
 \left( y\,x-1 \right) ^{2}}} \right) {M}+\\
 & \left( 1/2\,{\frac {
 \left( {z}^{2}+1 \right) ^{2}}{ \left( z\,x+1 \right) ^{2} \left( x-z
 \right) ^{2}}}+1/2\,{\frac { \left( {z}^{2}+1 \right) ^{2}}{ \left( z
+x \right) ^{2} \left( z\,x-1 \right) ^{2}}} \right) {m},\\
\dfrac{(y^2-1)y}{(1+y^2)^4}&-\dfrac{2yx(x^2+1)^2(y^2-1)(1-x^2)}{(x^2y^2-1)^2(x^2-y^2)^2}=-1/2\,{\frac {{\mu}}{{y}^{2}}}+1/8\,{\frac {{M}\, \left( {y}^{
2}+1 \right) ^{2}}{ \left( {y}^{2}-1 \right) ^{2}{y}^{2}}}\\
&+ \left( 1/2
\,{\frac { \left( {z}^{2}+1 \right) ^{2}}{ \left( z\,y+1 \right) ^{2}
 \left( -z+y \right) ^{2}}}+1/2\,{\frac { \left( {z}^{2}+1 \right) ^{2
}}{ \left( z+y \right) ^{2} \left( z\,y-1 \right) ^{2}}} \right) {m},\\
\dfrac{(z^2-1)z}{(1+z^2)^4}-&\dfrac{2zx(x^2+1)^2(1-x^2)(z^2-1)}{(x^2z^2-1)^2(x^2-z^2)^2}=-1/2\,{\frac {{\mu}}{{z}^{2}}}\\
& -\dfrac{2zy(y^2+1)^2(y^2-1)(z^2-1)}{(y^2z^2-1)^2(y^2-z^2)^2}{M}+1/8\,{\frac {{m}\,
 \left( {z}^{2}+1 \right) ^{2}}{ \left( {z}^{2}-1 \right) ^{2}{z}^{2}}}
\end{split}
\end{equation}

Let be
\begin{equation}
\begin{split}
A_1=&-\dfrac{(1-x^2)x}{(1+x^2)^4}+1/8\,{\frac { \left( {x}^{2}+1 \right) ^{2}}{ \left( {x}^{2}-1
 \right) ^{2}{x}^{2}}},\\
A_2=&\dfrac{(y^2-1)y}{(1+y^2)^4}-\dfrac{2yx(x^2+1)^2(y^2-1)(1-x^2)}{(x^2y^2-1)^2(x^2-y^2)^2},\\
A_3=&\dfrac{(z^2-1)z}{(1+z^2)^4}-\dfrac{2zx(x^2+1)^2(1-x^2)(z^2-1)}{(x^2z^2-1)^2(x^2-z^2)^2},\\
a_1=&-1/2\,{\frac {{1}}{{x}^{2}}},\\
a_2=&-1/2\,{\frac {{1}}{{y}^{2}}},\\
a_3=&-1/2\,{\frac {{1}}{{z}^{2}}},\\
b_1=&\left( 1/2\,{\frac { \left( {y}^{2}
+1 \right) ^{2}}{ \left( y\,x+1 \right) ^{2} \left( x-y \right) ^{2}}}
+1/2\,{\frac { \left( {y}^{2}+1 \right) ^{2}}{ \left( y+x \right) ^{2}
 \left( y\,x-1 \right) ^{2}}} \right),\\
b_2=&1/8\,{\frac { \left( {y}^{
2}+1 \right) ^{2}}{ \left( {y}^{2}-1 \right) ^{2}{y}^{2}}},\\
b_3=&-\dfrac{2zy(y^2+1)^2(y^2-1)(z^2-1)}{(y^2z^2-1)^2(y^2-z^2)^2},\\
c_1=&\left( 1/2\,{\frac {
 \left( {z}^{2}+1 \right) ^{2}}{ \left( z\,x+1 \right) ^{2} \left( x-z
 \right) ^{2}}}+1/2\,{\frac { \left( {z}^{2}+1 \right) ^{2}}{ \left( z
+x \right) ^{2} \left( z\,x-1 \right) ^{2}}} \right), \\
c_2=&\left( 1/2
\,{\frac { \left( {z}^{2}+1 \right) ^{2}}{ \left( z\,y+1 \right) ^{2}
 \left( -z+y \right) ^{2}}}+1/2\,{\frac { \left( {z}^{2}+1 \right) ^{2
}}{ \left( z+y \right) ^{2} \left( z\,y-1 \right) ^{2}}} \right), \\
c_3=&1/8\,{\frac {
 \left( {z}^{2}+1 \right) ^{2}}{ \left( {z}^{2}-1 \right) ^{2}{z}^{2}}}.\\
\end{split}
\end{equation}

We have $a_1,a_2,a_3,b_3<0$ and $A_1,b_1,b_2,c_1,c_2,c_3>0$. System (\ref{ss3}) becomes

\begin{equation}\label{syst5}
\begin{split}
A_1=a_1\mu+b_1M+c_1m,\\
A_2=a_2\mu+b_2M+c_2m,\\
A_3=a_3\mu+b_3M+c_3m.\\
\end{split}
\end{equation}

From the second equation, we have $\mu=\dfrac{A_2}{a_2}-\dfrac{b_2M}{a_2}-\dfrac{c_2m}{a_2}$. Substituting into the other equations we have

\begin{equation}
\begin{split}
a_2A_1-a_1A_2=M(b_1a_2-a_1b_2)+m(c_1a_2-a_1c_2),\\
a_2A_3-a_3A_2=M(b_3a_2-a_3b_2)+m(a_2c_3-a_3c_2).
\end{split}
\end{equation}

Adding these two equations,

\begin{equation}
\begin{split}
-a_2(A_1+A_3)+A_2(a_1+a_3)+M[a_2(b_1+b_3)-b_2(a_1+a_3)]\\
+m[a_2(c_1+c_3)-c_2(a_1+a_3)]=0. \label{Mmm}
\end{split}
\end{equation}

As in the previous case, among the functions $A_i,a_i,b_i,c_i, i=1,2,3$, the only values that depend on $x$ and $y$ simultaneously are $A_2$ and $b_1$. Fix $x$ and $z$, and consider values of $y$ close enough to $\frac{1}{x}$, in such a way that  $-a_2(A_1+A_3)+A_2(a_1+a_3)>0$ (since $A_2(a_1+a_3)>0$ ) and $a_2(b_1+b_3)-b_2(a_1+a_3)<0$ (since $a_2(b_1+b_3)<0$). This is enough to conclude the existence of $m,M>0$ such that  (\ref{Mmm}) holds. Notice that for this values of $x,y,z$ we have $\mu>0$.

\begin{itemize}
\item Case: $m_2, m_3,m_4,m_5$ lie inside the geodesic circle and $ m_6,m_7$ lie outside the geodesic circle, with $x<y<1<z<1/y$.
\end{itemize}
The equation corresponding to particle $q_4$ is

\begin{equation}\label{ee9}
\begin{split}
-\dfrac{(1-y^2)y}{(1+y^2)^4}&+1/2\,{\frac { \left( {x}^{2}+1 \right) ^{2}}{ \left( y\,x+1 \right) ^{
2} \left( -y+x \right) ^{2}}}+1/2\,{\frac { \left( {x}^{2}+1 \right) ^
{2}}{ \left( y+x \right) ^{2} \left( y\,x-1 \right) ^{2}}}
=-1/2\,{\frac {{\mu}}{{y}^{2}}}\\
&-1/8\,{\frac {M \left( {y}^{
2}+1 \right) ^{2}}{ \left( {y}^{2}-1 \right) ^{2}{y}^{2}}}-\dfrac{2yz(z^2+1)^2(z^2-1)(1-y^2)} {(y^2z^2-1)^2(y^2-z^2)^2} {
m},
\end{split}
\end{equation}

Left part of latter equation is positive by Lemma \ref{lema2}.  Hence equation (\ref{ee9})  is never satisfied for any $\mu,M,m>0$.

\begin{itemize}
\item Case: $m_2, m_3,m_4,m_5$ lie inside the geodesic circle and $ m_6,m_7$ lie outside the geodesic circle, with $x<y<1<1/y<z<1/x$
\end{itemize}

\begin{equation}\label{ss4}
\begin{split}
-\dfrac{(1-x^2)x}{(1+x^2)^4}+&1/8\,{\frac { \left( {x}^{2}+1 \right) ^{2}}{ \left( {x}^{2}-1
 \right) ^{2}{x}^{2}}}=-\dfrac{1}{2}\dfrac{\mu}{x^2}+\dfrac{2xy(y^2+1)^2(1-y^2)(1-x^2)}{(x^2y^2-1)^2(x^2-y^2)^2}M\\&-\dfrac{2xz(z^2+1)^2(z^2-1)(1-x^2)}{(x^2y^2-1)^2(x^2-z^2)^2}m,
\\
-\dfrac{(1-y^2)y}{(1+y^2)^4}&+1/2\,{\frac { \left( {x}^{2}+1 \right) ^{2}}{ \left( x\,y+1 \right) ^{
2} \left( -y+x \right) ^{2}}}+1/2\,{\frac { \left( {x}^{2}+1 \right) ^
{2}}{ \left( y+x \right) ^{2} \left( x\,y-1 \right) ^{2}}}
=-1/2\,{\frac {{\mu}}{{y}^{2}}}\\
&-1/8\,{\frac {M \left( {y}^{
2}+1 \right) ^{2}}{ \left( {y}^{2}-1 \right) ^{2}{y}^{2}}}+ \left( 1/2
\,{\frac { \left( {z}^{2}+1 \right) ^{2}}{ \left( z\,y+1 \right) ^{2}
 \left( y-z \right) ^{2}}}+1/2\,{\frac { \left( {z}^{2}+1 \right) ^{2}
}{ \left( z+y \right) ^{2} \left( z\,y-1 \right) ^{2}}} \right) {m},\\
\dfrac{(z^2-1)z}{(1+z^2)^4}&+1/2\,{\frac { \left( {x}^{2}+1 \right) ^{2}}{ \left( z\,x+1 \right) ^{
2} \left( x-z \right) ^{2}}}+1/2\,{\frac { \left( {x}^{2}+1 \right) ^{
2}}{ \left( z+x \right) ^{2} \left( z\,x-1 \right) ^{2}}}=-1/2\,{\frac {\mu}{{z}^{2}}}\\
&+2\,{\frac { \left( {y}^{2}+1
 \right) ^{2}z\,y\, \left( {z}^{2}-1 \right)  \left( -{y}^{2}+1
 \right) M}{ \left( {y}^{2}{z}^{2}+1 \right) ^{2} \left( {y}^{2
}-{z}^{2} \right) ^{2}}}+1/8\,{\frac {m \left( {z}^{2}+1
 \right) ^{2}}{ \left( {z}^{2}-1 \right) ^{2}{z}^{2}}}
.
\end{split}
\end{equation}

Let be

\begin{equation}
\begin{split}
A_1=&-\dfrac{(1-x^2)x}{(1+x^2)^4}+1/8\,{\frac { \left( {x}^{2}+1 \right) ^{2}}{ \left( {x}^{2}-1
 \right) ^{2}{x}^{2}}},\\
A_2=&-\dfrac{(1-y^2)y}{(1+y^2)^4}+1/2\,{\frac { \left( {x}^{2}+1 \right) ^{2}}{ \left( x\,y+1 \right) ^{
2} \left( -y+x \right) ^{2}}}+1/2\,{\frac { \left( {x}^{2}+1 \right) ^
{2}}{ \left( y+x \right) ^{2} \left( x\,y-1 \right) ^{2}}},\\
A_3=&\dfrac{(z^2-1)z}{(1+z^2)^4}+1/2\,{\frac { \left( {x}^{2}+1 \right) ^{2}}{ \left( z\,x+1 \right) ^{
2} \left( x-z \right) ^{2}}}+1/2\,{\frac { \left( {x}^{2}+1 \right) ^{
2}}{ \left( z+x \right) ^{2} \left( z\,x-1 \right) ^{2}}},\\
a_1=&-\dfrac{1}{2x^2},\\
a_2=&-\dfrac{1}{2y^2},\\
a_3=&-\dfrac{1}{2z^2},\\
b_1=&\dfrac{2xy(y^2+1)^2(1-y^2)(1-x^2)}{(x^2y^2-1)^2(x^2-y^2)^2},\\
b_2=&-1/8\,{\frac { \left( {y}^{
2}+1 \right) ^{2}}{ \left( {y}^{2}-1 \right) ^{2}{y}^{2}}},\\
b_3=&{\frac { 2zy\left( {y}^{2}+1
 \right) ^{2} \left( {z}^{2}-1 \right)  \left( 1-{y}^{2}
 \right)}{ \left( {y}^{2}{z}^{2}+1 \right) ^{2} \left( {y}^{2
}-{z}^{2} \right) ^{2}}},\\
c_1=&-\dfrac{2xz(z^2+1)^2(z^2-1)(1-x^2)}{(x^2z^2-1)^2(x^2-z^2)^2},\\
c_2=&\left( 1/2
\,{\frac { \left( {z}^{2}+1 \right) ^{2}}{ \left( z\,y+1 \right) ^{2}
 \left( y-z \right) ^{2}}}+1/2\,{\frac { \left( {z}^{2}+1 \right) ^{2}
}{ \left( z+y \right) ^{2} \left( z\,y-1 \right) ^{2}}} \right), \\
c_3=&1/8\,{\frac {\left( {z}^{2}+1
 \right) ^{2}}{ \left( {z}^{2}-1 \right) ^{2}{z}^{2}}}.\\
\end{split}
\end{equation}

The signs of $A_1$ and $A_2$ are given by lemmas (\ref{lemma5}) and (\ref{lema2}), respectively.

We have $A_1,A_2,A_3,b_1,b_2,c_2,c_3>0$; $a_1,a,_2,a_3,b_2,c_1<0$. Then system (\ref{ss4})  becomes

\begin{equation}
\begin{split}
A_1=a_1\mu+b_1M+c_1m,\\
A_2=a_2\mu+b_2M+c_2m,\\
A_3=a_3\mu+b_3M+c_3m.\\
\end{split}
\end{equation}

From the second equation we have $m=\frac{A_2}{c_2}-\frac{a_2}{c_2}\mu-\frac{b_2}{c_2}M.$ Substituting into the other two equations,

\begin{equation}
\begin{split}
c_2A_1-A_2c_1=\mu(a_1c_2-c_1a_2)+M(b_1c_2-c_1b_2),\\
c_2A_3-A_2c_3=\mu(a_3c_2-c_3a_2)+M(b_3c_2-c_3b_2).
\end{split}
\end{equation}

Adding the two equations,

\begin{equation} \label{Mmm2}
-c_2(A_1+A_3)+A_2(c_1+c_3)+\mu[c_2(a_1+a_3)-a_2(c_1+c_3)]+M[c_2(b_1+b_3)-b_2(c_1+c_3)]=0.
 \end{equation}

Notice that $c_2$ and $b_3$ are the only equations that depend simultaneously on $z$ and $y$. Fix $y$ and $x$. We can take $z$ in such a way that $c_2$ and $b_3$ are as large as we want, in this way we pick $z$ such that $-c_2(A_1+A_3)+A_2(c_1+c_3)<0$ and $c_2(b_1+b_3)-b_2(c_1+c_3)>0$. Then there exist $\mu$ and $M$ positive such that (\ref{Mmm2}) holds.

\begin{itemize}
\item Case: $m_2, m_3,m_4,m_5$ lie inside the geodesic circle and $ m_6,m_7$ lie outside the geodesic circle, with $x<y<1<1/x<z$
\end{itemize}
\begin{equation}\label{ss6}
\begin{split}
-\dfrac{(1-x^2)x}{(1+x^2)^4}+&1/8\,{\frac { \left( {x}^{2}+1 \right) ^{2}}{ \left( {x}^{2}-1
 \right) ^{2}{x}^{2}}}=-1/2\,{\frac {\mu}{{x}^{2}}}+2\,{\frac { \left( {y}^{2}+1 \right) ^{2}x\,y\, \left( -{y}^{2}+1
 \right)  \left( -{x}^{2}+1 \right) M}{ \left( {x}^{2}-{y}^{2}
 \right) ^{2} \left( {x}^{2}{y}^{2}-1 \right) ^{2}}}\\
 &+\left( 1/2\,{\frac { \left( {z}^{2}+1 \right) ^{2}}{ \left( z\,x+1
 \right) ^{2} \left( x-z \right) ^{2}}}+1/2\,{\frac { \left( {z}^{2}+1
 \right) ^{2}}{ \left( z+x \right) ^{2} \left( z\,x-1 \right) ^{2}}}
 \right) m,\\
-\dfrac{(1-y^2)y}{(1+y^2)^4}+&1/2\,{\frac { \left( {x}^{2}+1 \right) ^{2}}{ \left( x\,y+1 \right) ^{
2} \left( -y+x \right) ^{2}}}+1/2\,{\frac { \left( {x}^{2}+1 \right) ^
{2}}{ \left( y+x \right) ^{2} \left( x\,y-1 \right) ^{2}}}=-1/2\,{\frac {\mu}{{y}^{2}}}
\\
&-1/8\,{\frac {M \left( {y}^{2}+1
 \right) ^{2}}{ \left( {y}^{2}-1 \right) ^{2}{y}^{2}}}+ \left( 1/2\,{
\frac { \left( {z}^{2}+1 \right) ^{2}}{ \left( z\,y+1 \right) ^{2}
 \left( y-z \right) ^{2}}}+1/2\,{\frac { \left( {z}^{2}+1 \right) ^{2}
}{ \left( z+y \right) ^{2} \left( z\,y-1 \right) ^{2}}} \right) m,\\
\dfrac{(z^2-1)z}{(1+z^2)^4}-&2\,{\frac {x\,z\, \left( {x}^{2}+1 \right) ^{2} \left( {z}^{2}-1
 \right)  \left( -{x}^{2}+1 \right) }{ \left( {x}^{2}-{z}^{2} \right) 
^{2} \left( {x}^{2}{z}^{2}-1 \right) ^{2}}}=\\
&-1/8\,{\frac {\mu}{{z}^{2}}}+2\,{\frac {y\,z\, \left( {y}^{2}+1
 \right) ^{2} \left( {z}^{2}-1 \right)  \left( -{y}^{2}+1 \right) M}{
 \left( {y}^{2}-{z}^{2} \right) ^{2} \left( {y}^{2}{z}^{2}-1 \right) ^
{2}}}+1/8\,{\frac {m \left( {z}^{2}+1 \right) ^{2}}{ \left( {z}^{2}-1
 \right) ^{2}{z}^{2}}}.
\end{split}
\end{equation}

Let be

\begin{equation}
\begin{split}
A_1&=-\dfrac{(1-x^2)x}{(1+x^2)^4}+1/8\,{\frac { \left( {x}^{2}+1 \right) ^{2}}{ \left( {x}^{2}-1
 \right) ^{2}{x}^{2}}},\\
A_2&=-\dfrac{(1-y^2)y}{(1+y^2)^4}+1/2\,{\frac { \left( {x}^{2}+1 \right) ^{2}}{ \left( x\,y+1 \right) ^{
2} \left( -y+x \right) ^{2}}}+1/2\,{\frac { \left( {x}^{2}+1 \right) ^
{2}}{ \left( y+x \right) ^{2} \left( x\,y-1 \right) ^{2}}},\\
A_3&=\dfrac{(z^2-1)z}{(1+z^2)^4}-2\,{\frac {x\,z\, \left( {x}^{2}+1 \right) ^{2} \left( {z}^{2}-1
 \right)  \left( -{x}^{2}+1 \right) }{ \left( {x}^{2}-{z}^{2} \right) 
^{2} \left( {x}^{2}{z}^{2}-1 \right) ^{2}}},\\
a_1&=-1/2\,{\frac {1}{{x}^{2}}},\\
a_2&=-1/2\,{\frac {1}{{y}^{2}}},\\
a_3&=-1/2\,{\frac {1}{{z}^{2}}},\\
b_1&=2\,{\frac { \left( {y}^{2}+1 \right) ^{2}x\,y\, \left( -{y}^{2}+1
 \right)  \left( -{x}^{2}+1 \right) }{ \left( {x}^{2}-{y}^{2}
 \right) ^{2} \left( {x}^{2}{y}^{2}-1 \right) ^{2}}},\\
b_2&=-1/8\,{\frac { \left( {y}^{2}+1
 \right) ^{2}}{ \left( {y}^{2}-1 \right) ^{2}{y}^{2}}},\\
b_3&=2\,{\frac {y\,z\, \left( {y}^{2}+1
 \right) ^{2} \left( {z}^{2}-1 \right)  \left( -{y}^{2}+1 \right) }{
 \left( {y}^{2}-{z}^{2} \right) ^{2} \left( {y}^{2}{z}^{2}-1 \right) ^
{2}}},\\
c_1&=\left( 1/2\,{\frac { \left( {z}^{2}+1 \right) ^{2}}{ \left( z\,x+1
 \right) ^{2} \left( x-z \right) ^{2}}}+1/2\,{\frac { \left( {z}^{2}+1
 \right) ^{2}}{ \left( z+x \right) ^{2} \left( z\,x-1 \right) ^{2}}}
 \right),\\
c_2&=\left( 1/2\,{
\frac { \left( {z}^{2}+1 \right) ^{2}}{ \left( z\,y+1 \right) ^{2}
 \left( y-z \right) ^{2}}}+1/2\,{\frac { \left( {z}^{2}+1 \right) ^{2}
}{ \left( z+y \right) ^{2} \left( z\,y-1 \right) ^{2}}} \right),\\
c_3&=1/8\,{\frac { \left( {z}^{2}+1 \right) ^{2}}{ \left( {z}^{2}-1
 \right) ^{2}{z}^{2}}}.\\
\end{split}
\end{equation}

The signs of $A_1$ and $A_2$ are given by lemmas (\ref{lemma5}) and (\ref{lema2}), respectively.

We have $A_1,A_2,b_1,b_3,c_1,c_2,c_3>0$ and $,a_1,a_2,a_3,b_2<0$. The system (\ref{ss6}) becomes

\begin{equation}
\begin{split}
A_1=a_1\mu+b_1M+c_1m,\\
A_2=a_2\mu+b_2M+c_2m,\\
A_3=a_3\mu+b_3M+c_3m.\\
\end{split}
\end{equation}

From the second equation we have $m=\dfrac{A_2}{c_2}-\dfrac{a_2\mu}{c_2}-\dfrac{b_2M}{c_2}>0$. Substituting into the other two equations we have

\begin{equation}
\begin{split}
c_2A_1-c_1A_2=\mu(a_1c_2-c_1a_2)+M(b_1c_2-c_1b_2),\\
c_2A_3-c_3A_2=\mu(a_3c_2-c_3a_2)+M(b_3c_2-c_3b_2).
\end{split}
\end{equation}

Adding last two equations,

\begin{equation}\label{Mm4}
\begin{split}
-c_2(A_1+A_3)+A_2(c_1+c_3)&+\mu[c_2(a_1+a_3)-a_2(c_1+c_3)]\\
&+M[c_2(b_1+b_3)-b_2(c_1+c_3)]=0.
\end{split}
\end{equation}

Notice that $c_2(b_1+b_3)-b_2(c_1+c_3)>0$. If we take $1/z$ close enough to 0, or $z$ large enough, we have that $a_3$ is very large in absolute value, hence in this way we can choose $z$ such that $c_2(a_1+a_3)-a_2(c_1+c_3)<0$. At this point we can conclude that there exist $\mu$ and $M$ positive such that (\ref{Mm4}) holds.

With all the above we have finished the proof of Theorem \ref{Theorem 2}.

\bigskip

\begin{cor}
In the 6-body problem on $\mathbb{M}^2$ we consider 6 particles on the same geodesic with masses $m_1=m_2=1$, $m_3=m_4=M$ and $m_5=m_6=m$, in a symmetric configuration  with initial positions $z_1=-z_2=x>0$, $z_3=-z_4=y>0$ and $z_5=-z_6=z>0$ ($x<y<z$).
\begin{itemize}
\item If $m_2,m_3, m_4,m_5,m_6,m_7$ lie inside the geodesic circle, then do not exist relative equilibria.
\item If $m_2,m_3$ lie inside the geodesic circle, and $ m_4,m_5,m_6,m_7$ lie outside the geodesic circle with $y<z<\frac{1}{x}$, then do not exist relative equilibria. If $y<1/x<z$ or $1/x<y<z$ then it is possible to find relative equilibria.
\end{itemize}
\end{cor}

\begin{proof}

The proof is similar as in Lemma \ref{4b}. Consider $\mu=0$.\\

\begin{itemize}
\item If $m_2,m_3, m_4,m_5,m_6,m_7$ lie inside the geodesic circle, then this case correspond to $x<y<z<1$.
\end{itemize}

 Using condition (\ref{cond}) for particle $z_6$ we get,

\begin{equation}
\begin{split}
-\dfrac{(1-z^2)z}{(1+z^2)^4} & +\dfrac{1}{2}(x^2+1)^2\left(\dfrac{1}{(xz+1)^2(x-z)^2}+\dfrac{1}{(z+x)^2(xz-1)^2}  \right)=\\
&-\dfrac{1}{2}(y^2+1)^2\left( \dfrac{1}{(yz+1)^2}+\dfrac{1}{(z+y)^2(yz-1)^2} \right)M-\dfrac{1}{8}\dfrac{(z^2+1)^2m}{(z^2-1)^2z^2}.
\end{split}
\end{equation}

Last equation is never satisfied for $M,m>0$, since left part is positive.

\begin{itemize}
\item Now we analyze the case when $m_2, m_3$ lie inside the geodesic circle and$ m_4,m_5,m_6,m_7$ lie outside the geodesic circle, with $y<z<1/x$.
\end{itemize}
The equation for this case corresponding to particle $m_2$ is

\begin{equation}\label{rr1}
\begin{split}
-\dfrac{(1-x^2)x}{(1+x^2)^4}&+1/8\,{
\frac { \left( {x}^{2}+1 \right) ^{2}}{ \left( {x}^{2}-1 \right) ^{2}{
x}^{2}}}=\\
&+ \left( 1/2\,{\frac { \left( {y}^{2}
+1 \right) ^{2}}{ \left( xy+1 \right) ^{2} \left( -y+x \right) ^{2}}}-
1/2\,{\frac { \left( {y}^{2}+1 \right) ^{2}}{ \left( y+x \right) ^{2}
 \left( xy-1 \right) ^{2}}} \right) {M}\\
 &+ \left( 1/2\,{\frac {
 \left( {z}^{2}+1 \right) ^{2}}{ \left( xz+1 \right) ^{2} \left( x-z
 \right) ^{2}}}-1/2\,{\frac { \left( {z}^{2}+1 \right) ^{2}}{ \left( z
+x \right) ^{2} \left( xz-1 \right) ^{2}}} \right) {m}.\\
\end{split}
\end{equation}

The factors for $m$ and $M$ in the last equation can be seen as

\[ \dfrac{1}{2}\dfrac{yx(y^2+1)^2(y^2-1)(x^2-1)}{(xy+1)^2(x^2-y^2)^2(xy-1)^2}<0, \]
\[  \dfrac{1}{2}\dfrac{zx(z^2+1)^2(z^2-1)(x^2-1)}{(xz+1)^2(x^2-z^2)^2(xz-1)^2}<0.
\]

Hence (\ref{rr1}) is never satisfied, since left part of the equation is positive.

\begin{itemize}
\item Case $y<1/x<z$.
\end{itemize}

Considering $\mu=0$, system (\ref{system1}) becomes

\begin{equation}
\begin{split}
A_1=b_1M+c_1m,\\
A_2=b_2M+c_2m,\\
A_3=b_3M+c_3m.\\
\end{split}
\end{equation}

From second equation $m=\frac{A_1}{c1}-\frac{b_1M}{c1}>0$. Substituting into the other two equations and adding them we have

\begin{equation}\label{s7}
\begin{split}
-c_1(A_2+A_3)+A_1(c_2+c_3)+M[c_1(b_2+b_3)-b_1(c_2+c_3)]=0.
\end{split}
\end{equation}

Among the functions $A_i,a_i,b_i,c_i, i=1,2,3$, the only values that depend on $x$ and $y$ simultaneously are $A_2$ and $b_1$. Consider values of $y$ close enough to $\frac{1}{x}$, in such a way that $-c_1(A_2+A_3)+A_1(c_2+c_3)<0$ and $[c_1(b_2+b_3)-b_1(c_2+c_3)]>0$ are satisfied.  When these two inequalities are fulfilled, then we can conclude that there are values for $M$ such that equation (\ref{s7}) is valid.

\begin{itemize}
\item Let us analyze the case when $m_2, m_3$ lie inside the geodesic circle and$ m_4,m_5,m_6,m_7$ lie outside the geodesic circle, with $1/x<y<z$.
\end{itemize}

System (\ref{syst5}) becomes

\begin{equation}
\begin{split}
A_1=b_1M+c_1m,\\
A_2=b_2M+c_2m,\\
A_3=b_3M+c_3m.\\
\end{split}
\end{equation}

With $b_3<0$, and $A_1,b_1,b_2,c_1,c_2,c_3>0$. It is important to notice that $A_2$ and $A_3$ might be positive or negative, and last system makes sense only if $A_2,A_3>0$. For $y,z$ fixed, consider $x$ close enough to $0$ in such a way that $A_2,A_3>0$, for this value of $x$, $A_1$ is very large.

 From third equation from the above system $\frac{A_3}{c_3}-\frac{b_3M}{c_3}=m>0$. Substituting into the other two equations and adding them we have

\begin{equation}
c_3(A_1+A_2)-A_3(c_1+c_2)=M[c_3(b_1+b_2)-b_3(c_1+c_2)].
\end{equation}

It is easy to check that right part of last equation is positive. Left part is positive because we have chosen $x$ such that $A_1$ is large enough. Hence, there exist $M>0$ such that last equation holds.

\end{proof}

We can generalize the above result for $n$ masses with symmetric configuration for the case  where no relative equilibria exist.

\begin{prop}
Consider $n$ (odd) particles on $\mathbb{M}^2$. We consider particles on the same geodesic with masses $m_1$, $m_2=m_3=1, \ m_4=m_5, \dots, m_{n-1}=m_n$ and initial positions $0=z_1<z_2=-z_3<z_4=-z_5<\dots<z_{n-1}=-z_n.$ Then do not exist relative equilibria if
\begin{itemize}
\item The $n$ particles are inside the geodesic circle.
\item All the particles except the bodies 2 and 3 are outside the geodesic circle with $z_{n-1}<1/z_2$
\item All the particles except the bodies $n-1$ and $n$ are inside the geodesic circle with $z_{n-1}<1/z_{n-3}$.
\end{itemize}

\end{prop}

\begin{proof}
\begin{itemize}
\item Case $z_{n-1}<1$. 
\end{itemize}
Equation for particle $n-1$ in condition (\ref{cond}) is

\begin{equation}\label{exx}
\begin{split}
-\dfrac{(1-z_{n-1}^2)z_{n-1}}{(1+z_{n-1}^2)^4}=&-\sum_{i=1,i\neq n-1}^n\dfrac{m_i(z_i^2+1)^2}{2(1+z_iz_{n-1})^2(z_i-z_{n-1})^2}\\
=&-\sum_{i=1, \ i\neq n-1, \ i\neq 2,3}^n\dfrac{m_i(z_i^2+1)^2}{2(1+z_iz_{n-1})^2(z_i-z_{n-1})^2}\\
&-\dfrac{(z_2^2+1)^2}{2(1+z_2z_{n-1})^2(z_2-z_{n-1})^2}-\dfrac{(z_3^2+1)^2}{2(1+z_3x_{n-1})^2(z_3-z_{n-1})^2}.
\end{split}
\end{equation}

By Lemma (\ref{lema2}) we have 

\[ -\dfrac{(1-z_{n-1}^2)z_{n-1}}{(1+z_{n-1}^2)^4}+\dfrac{(z_2^2+1)^2}{2(1+z_2z_{n-1})^2(z_2-z_{n-1})^2}+\dfrac{(z_3^2+1)^2}{2(1+z_3x_{n-1})^2(z_3-z_{n-1})^2} >0.\]

Hence expression (\ref{exx}) has no solution for $m_i>0$.

\begin{itemize}
\item Case $z_{n-1}>z_4>1>z_2$ and $z_n<1/z_2$.
\end{itemize}

Equation for particle $2$ in condition (\ref{cond}) is

\begin{equation}
\begin{split}
-\dfrac{(1-z_{2}^2)z_{2}}{(1+z_{2}^2)^4}=&\sum_{i=1,i\neq 2}^n\dfrac{m_i(z_i^2+1)^2(1+z_iz_{2})(z_i-z_{2})}{2|(1+z_iz_{2})|^3|(z_i-z_{2})|^3}\\
=&\sum_{i=1,i\neq 2, 3}^n\dfrac{m_i(z_i^2+1)^2(1+z_iz_{2})(z_i-z_{2})}{2|(1+z_iz_{2})|^3|(z_i-z_{2})|^3}-\dfrac{1}{8}\dfrac{(z_2^2+1)^2}{(1-z_2^2)^2z_2^2}.
\end{split}
\end{equation}

We can write last equation as

\begin{equation}
\begin{split}
-\dfrac{(1-z_{2}^2)z_{2}}{(1+z_{2}^2)^4}+ \dfrac{1}{8}\dfrac{(z_2^2+1)^2}{(1-z_2^2)^2z_2^2}=&\sum_{i=1,i\neq 2, 3}^n\dfrac{m_i(z_i^2+1)^2(1+z_iz_{2})(z_i-z_{2})}{2|(1+z_iz_{2})|^3|(z_i-z_{2})|^3}.
\end{split}
\end{equation}

Left part is positive by lemma (\ref{lemma5}). Right part can be seen as

\begin{equation}\label{exp5}
\begin{split}
 &-\dfrac{1/2m_1}{z_2}+\dfrac{m_4}{2}\left[\dfrac{(z_4^2+1)^2}{(1+z_4z_{2})^2(z_4-z_{2})^2} -\dfrac{(z_4^2+1)^2}{(z_4z_{2}-1)^2(z_4+z_{2})^2} \right]\\
 &+\dfrac{m_6}{2}\left[\dfrac{(z_6^2+1)^2}{(1+z_6z_{2})^2(z_6-z_{2})^2} -\dfrac{(z_6^2+1)^2}{(z_6z_{2}-1)^2(z_6+z_{2})^2} \right] + \cdots\\
&+\dfrac{m_{n-1}}{2}\left[\dfrac{(z_{n-1}^2+1)^2}{(1+z_{n-1}z_{2})^2(z_{n-1}-z_{2})^2} -\dfrac{(z_{n-1}^2+1)^2}{(z_{n-1}z_{2}-1)^2(z_{n-1}+z_{2})^2} \right].
 \end{split}
 \end{equation}

We have 

\begin{equation}\label{negative}
 \left[\dfrac{(z_i^2+1)^2}{(1+z_iz_{2})^2(z_i-z_{2})^2} -\dfrac{(z_i^2+1)^2}{(z_iz_{2}-1)^2(z_i+z_{2})^2} \right]=\dfrac{z_iz_2(z_i^2+1)^2(z_i^2-1)(z_2-1)}{(z_iz_2+1)^2(z_i-z^2)^2(z_1z_2-1)^2}<0,
 \end{equation}
for $i=4,6,8,\dots,n-1$. Hence expression (\ref{exp5}) has no solution for $m_i>0$.

\begin{itemize}
\item Case $z_{n-1}>1>z_{n-3}>z_2$ and $z_n<1/z_{n-3}$.
\end{itemize}

Equation for particle $n-3$ in condition (\ref{cond}) is

\begin{equation}
\begin{split}
-\dfrac{(1-z_{n-3}^2)z_{n-3}}{(1+z_{n-3}^2)^4}=&\sum_{i=1,i\neq n-3}^n\dfrac{m_i(z_i^2+1)^2(1+z_iz_{2})(z_i-z_{2})}{2|(1+z_iz_{2})|^3|(z_i-z_{2})|^3}\\
=&-\sum_{i=1,i\neq n-1}^{n-2}\dfrac{m_i(z_i^2+1)^2}{2(1+z_iz_{n-1})^2(z_i-z_{n-1})^2}\\
&+\dfrac{m_{n-1}}{2}\left[\dfrac{(z_{n-1}^2+1)^2}{(1+z_{n-1}z_{2})^2(z_{n-1}-z_{2})^2} -\dfrac{(z_{n-1}^2+1)^2}{(z_{n-1}z_{2}-1)^2(z_{n-1}+z_{2})^2} \right]\\
=&-\sum_{i=1,i\neq 2,3,n-1}^{n-2}\dfrac{m_i(z_i^2+1)^2}{2(1+z_iz_{n-1})^2(z_i-z_{n-1})^2}\\
&-\dfrac{(z_2^2+1)^2}{2(1+z_2z_{n-1})^2(z_2-z_{n-1})^2}-\dfrac{(z_3^2+1)^2}{2(1+z_3x_{n-1})^2(z_3-z_{n-1})^2}\\
&+\dfrac{m_{n-1}}{2}\left[\dfrac{(z_{n-1}^2+1)^2}{(1+z_{n-1}z_{2})^2(z_{n-1}-z_{2})^2} -\dfrac{(z_{n-1}^2+1)^2}{(z_{n-1}z_{2}-1)^2(z_{n-1}+z_{2})^2} \right].
\end{split}
\end{equation}

Last expression can be seen as

\begin{equation}
\begin{split}
&-\dfrac{(1-z_{n-3}^2)z_{n-3}}{(1+z_{n-3}^2)^4}+\dfrac{(z_2^2+1)^2}{2(1+z_2z_{n-1})^2(z_2-z_{n-1})^2}+\dfrac{(z_3^2+1)^2}{2(1+z_3x_{n-1})^2(z_3-z_{n-1})^2}=\\
&-\sum_{i=1,i\neq 2,3,n-1}^{n-2}\dfrac{m_i(z_i^2+1)^2}{2(1+z_iz_{n-1})^2(z_i-z_{n-1})^2}\\
&+\dfrac{m_{n-1}}{2}\left[\dfrac{(z_{n-1}^2+1)^2}{(1+z_{n-1}z_{2})^2(z_{n-1}-z_{2})^2} -\dfrac{(z_{n-1}^2+1)^2}{(z_{n-1}z_{2}-1)^2(z_{n-1}+z_{2})^2} \right].
\end{split}
\end{equation}

In last equation, left part is positive (Lemma \ref{lema2}) and right part is negative (as seen in equation \ref{negative} ).
\end{proof}

\subsection*{Acknowledgements} The first author has been partially supported by {\it Asosiaci\'on Mexicana de Cultura A.C.}  The second author was supported by {\it The 2017's Plan of Foreign Cultural and Educational Experts Recruitment for the Universities Under the Direct Supervision of the Ministry of Education of China} (Grant no. WQ2017SCDX045). Part of this work was done during the visit to the first author to Sichuan University, Chengdu, China. We thank to professor S. Zhang for the invitation and hospitality.

 \end{document}